# Prime Ideals in Noetherian Rings

C. L. Wangneo

# A-1, Jammu university,Rail Road Jammu,Jammu & Kashmir,India(180006).

**Abstract**: In this short note we study the links of certain prime ideals of a noetherian ring R. We first give the definition of a link krull symmetric noetherian ring R. We then prove theorem 9 that states that for any linked prime ideals P ′ and Q′ of the polynomial ring R[X] where R is a link krull symmetric noetherian ring, if The prime ideal P ′ is extended then Q′ is also an extended prime ideal of R[X]. An application of theorem 9 is then given in theorem 12 for the ring R[X] when R is assumed to be a fully bounded noetherian ring.

**1. Introduction**: In a noetherian ring R there is a link from Q → P , for Q, P prime ideals of R if there is an ideal A of R Such that QP ≤ A < Q ∩ P and Q ∩ P/A is torsionfree as a left R/Q and as a right R/P module. In this paper we study the links between certain prime ideals of the polynomial ring R[X], where R is a link k-symmetric noetherian ring. We first state the definition of a link krull symmetric(link k-symmetric for short) noetherian ring R. We then prove the main theorem of this paper, namely theorem 9, which states that for the linked prime ideals Q and P of the the polynomial ring R[X], where R is a link k-symmetric noetherian ring, if the ideal P is an extended prime ideal of the ring R[X] then Q is also an extended prime ideal of the ring R[X]. Theorem 12 then is an application of theorem 9 to the case of a polynomial ring R[X], when R is considered to be a fully bounded noetherian ring. Moreover theorem 9 also provides easily in theorem 12 a description of all the links between any two prime ideals of the polynomial ring R[X] over a link k-symmetric noetherian ring.

**2. Definitions and Notation**: We mention that throughout all our rings
are with identity and all modules are unitary. A ring R is noetherian means that R is a right as well as a left noetherian ring. For the definition of right Krull dimension of a right R module M we refer the reader to [3 ] or [4]. A right module M is said to be right k-homogenous (k-homogenous for short) if every non zero sub module of M has same krull dimension as that of M . A semi- prime ideal S of a noetherian ring R is said to be right krull homogenous if the ring R/S is a right k-homogenous ring. A noetherian bimodule M over a ring R is said to be krull symmetric if the right krull dimension of M equals its left krull dimesion.For the definition of weak ideal invariance of an ideal A of a noetherian R , see[1]. For the definition of a fully bounded noetherian (FBN for short) ring, see [3]. We now mention a few words regarding our terminology in the present paper. For a ring R, Spec.R denotes the set of prime ideals of R and minSpec.R denotes the set of minimal prime ideals of R. For a krull symmetric bimodule M over a ring R, |M| denotes its left or right krull dimension. Finally if A is an ideal of a ring R then c(A) denotes the set of elements of R that are regular modulo A. For a ring R we denote by N (R) the nil radical of R.

**3. The Main Theorem**: We start with some Lemmas which help us in proving our main theorem.

Lemma 1. Let R be a Noetherian ring and let R[X] be the polynomial ring over R in a commuting indeterminate X. Let Q′ε spec.R[X] be a prime ideal such that Q′ ∩ R = Q. If c ε c(Q), then c ε c(Q′).

**Proof.** First observe that Q is a prime ideal of R. Next we consider two cases;

**Case 1** Q′ = Q[X]. In this case it is not difficult to see that c ε c(Q) implies that c ε c(Q′) also.

**Case 2** $Q' > Q[X]$. In this case consider the ring $R[X]/Q[X]$. Then $Q'/Q[X]$ is a nonzero prime ideal of the ring $R[X]$. If we prove that $c+Q[X]$, where $c \varepsilon c(Q)$, is an element such that $c \varepsilon c(Q')$, then the proof follows. Hence we may assume that $Q[X]=0$ and that R as well as $R[X]$ are prime rings. Thus $Q'$ is then a nonzero prime ideal of the ring $R[X]$ such that $Q' \cap R = 0$. In this context then c is a nonzero regular element of R. We now prove that $c \varepsilon c(Q')$. To see this, suppose $f(X) \varepsilon R[X]$ is an element such that $f(x)$ is not an element of $Q'$ and $cf(X) \varepsilon Q'$. Now it follows from Goldie [2], theorem (13) that since R is a Noetherian prime ring then the set $c(0)$ = set of regular elements of R, forms a left and a right ore-set in R. It is not difficult to see that $c(0)$ is a left and a right ore- sub set of the ring $R[X]$ too. Thus for any $g(X) \varepsilon R[X]$, we have that there exists an element $d \varepsilon c(0)$ such that $dg(X) \varepsilon R[X].c$ (after using the fact that $c(0)$ is a left ore-set of $R[x]$). Let $dg(X) = h(X).c$, for some $h(X) \varepsilon R[X]$. So we get that $dg(X)f(X) = h(X)cf(X)$, for all $g(X) \varepsilon R[X]$. Since it is given that $cf(X) \varepsilon Q'$, this thus implies that $dg(X)f(X)R[X] \leq Q'$, for all $g(X)\varepsilon R[X]$. This immediately implies that the $R[X]$ bi-module say,
$A = R[X]f[X]R[X] + Q'/Q'$ is left $c(0)$ torsion module. Since $R[X]$ is a Noethrian ring because R is a Noetherian ring, so the bi-module A is left and right finitely generated over $R[X]$. Since $R[X]$ is a prime ring, and A is $c(0)$ left torsion bi-module which is also finitely generated on either side, hence there exists an element $d_1 \varepsilon c(0)$ such that $d_1 A=0$. It is not difficult to see that this implies that $d_1 R[X]f[X] \varepsilon Q'$. Since $Q'$ is a prime ideal of $R[X]$, so either $d_1 \varepsilon Q'$ or $f(X) \varepsilon Q'$ which is not true. Thus $cf(X) \varepsilon Q'$ implies that $f(X) \varepsilon Q'$. Hence $c \varepsilon c(Q')$.

**Lemma 2.** Let R be a Noetherian ring and let $R[X]$ be the usual polynomial ring over R in a commuting indeterminate X. Let P be a prime ideal of R. Let $Q' \varepsilon$ spec.$R[X]$ and $B'$ be an ideal of $R[X]$ such that $Q'P[X] \leq B' < Q' \cap P[X]$. Let $Q = Q' \cap R$ and $B = B' \cap R$. Suppose further that $Q' \cap P[X]/B'$ is a nonzero torsionfree left $R[X]/Q'$ bi-module. Then the following hold:
(a) If $Q' \geq Q[X]$, then $Q \cap P/B$ is a left-$R/Q$ torsionfree module provided $Q \cap P/B$ is a nonzero R- bi-module.

(b) If $Q' = Q[X]$, then $Q \cap P/B$ is a nonzero R-bi-module that is left $R/Q$- torsionfree.

**Proof.** (a) First note that B is an ideal of R and Q is a prime ideal of R. We now prove (a) given that $Q \cap P/B \neq 0$. Suppose $Q \cap P/B$ is left $R/Q$ torsion module. Then there exists an element $a \varepsilon Q \cap P$ ( a is not an element of B) and an element $d \varepsilon c(Q)$, such that $da \varepsilon B$. Since a is not an element of B, so a is not an element of $B'$ also. Hence we get that there is an element $a \varepsilon Q' \cap P[X]$ and a is not an element of $B'$. By the previous lemma we have that $d \varepsilon c(Q')$. So $da \varepsilon B$ and $B \leq B'$ yields that $a + B'$ is a $c(Q)$ torsion element of the bi- module $Q' \cap P[X]/B'$. This ontradicts our hypothesis in (a). Hence $Q \cap P/B$ is a left-$R/Q$ torsionfree module.

(b) For the proof of (b) apply (a) above after observing that $Q' = Q[X]$ and $B' < Q' \cap P[X]$ implies immediately that $B < Q \cap P$ and hence $Q \cap P/B$ is a nonzero R bi-module.

**Lemma 3.** Let R be a Noetherian ring and let $R[X]$ be the usual Polynomial ring. Let $Q' \varepsilon$ spec.$R[X]$. Let $f(X) = a_0 + a_1 x_1 + \ldots + a_n X_n$ be a polynomial such that $f(X)$ is not an element of $Q'$. Assume $a_n$ is not an element of Q where $Q = Q' \cap R$. Then the following hold;
(a) If $f(X) \varepsilon c(Q')$, then there exists an element $r \varepsilon R$ such that $a_n r \varepsilon c(Q)$
(b) If $a_n \varepsilon c(Q)$, then $f(X) \varepsilon c(Q)[X]$ .
Proof. (a). The proof of (a) is on the same lines as that of Small [6,Lemma4.2].

(b). The proof of (b) is obvious.

**Lemma 4.** Let R be a Noetherian ring and let $R[X]$ be the usual Polynomial ring. Let Q and P be prime ideals of R and let B be an ideal of R such that $B < Q \cap P$ and $QP \leq B$. If $Q \cap P/B$ is a nonzero R-bi-module that is left $R/Q$ torsionfree module then $Q[X] \cap P[X]/B[X]$ is a nonzero $R[X]$-bi-module that is left $R[X]/Q[X]$ torsionfree module and conversely.

**Proof.** First note that since $Q \cap P/B$ is a nonzero R-bi-module hence

Q[X]∩P[X]/B[X] is also a nonzero R[X]-bi-module. Suppose Q[X]∩P[X]/B[X] is not a left R[X]/Q[X] torsionfree module. Then there exists an element f(X) in Q[X] ∩ P[X] (such that f(X) is not an element of B[X]) and an element d(X) of c((Q[X]) such that d(X)f(X) ε B[X]. Let a and d be the leading terms of f(X) and d(X) such that a and d do not belong to B and Q respectively. By the left version of Lemma (3) there exists an element r in R such that rd is an element of c(Q). Now d(X)f(X) ε B[X] implies that rda ε B. Note that a ε Q ∩ P and a is not an element of B. Since rd ε c(Q) we get thus that Q ∩ P/B is not a left R/Q torsionfree module. This is a contradiction to our hypothesis. Hence Q[X]∩P[X]/B[X] must be a left R[X]/Q[X] torsionfree module. The converse follows from Lemma 2.

**Definition 5**. Two prime ideals P and Q in a Noetherian ring R are said to have a second layer link, written Q → P, if there exists an ideal A of R with QP ≤ A and A < Q ∩ P such that Q ∩ P/A is a torsionfree R/Q − R/P bimodule. In this case we also say Q is linked to P or that Q is a second layer link to P. We also say in this case that Q ∩P/A is a linking bimodule for the link Q → P or Q is linked to P via the ideal A.

**Definition 6**. We say a Noetherian ring R is link krull symmetric (link k-symmetric for short) if for any prime ideals Q and P of R such that Q is linked to P we have |R/P| = |R/Q|.

**Theorem 7**. Let R be a Noetherian ring. Let R[X] be the polynomial ring over R. Then Q → P is a link between the prime ideals Q and P of R if and only if Q[X] → P[X] is a link between the corresponding prime ideals Q[X] and P[X] of R[X].

**Proof**. Let Q → P be a link between the prime ideals P and Q of R with Q∩P/B the corresponding linking bimodule. Then by Lemma (4) we get that Q[X] ∩ P[X]/B[X] is a linking bimodule for the prime ideals Q[X] and P[X] of R[X]. Hence Q[X] is linked to P[X]. Conversely let Q[X] → P[X] be a link of the prime ideals Q[X] and P[X] of R[X] via the bimodule Q[X]∩ P[X]/B′, where B′ is an ideal of R[X]. Suppose B = B′ ∩ R. It is clear that Q[X] → P[X] is also a link of the prime ideal Q[X] to P[X] via the via the ideal B[X]. Using Lemma 2 we get that Q→P is a link of the prime ideals Q and P via the ideal B.

**Theorem 8**. Let R be a link k-symmetric, Noetherian ring. Let Q → P be a link of the prime ideals Q and P of R via an ideal B of R. Then R/B has an artinian quotient ring.

**Proof**. Given that QP ≤ B < Q ∩ P. Also |R/P| = |R/Q| because R is a link k-symmetric ring. Thus P/B and Q/B are distinct incomparable prime ideals of R/B unless Q = P. Now it is also obvious that Q/B and P/B are the minimal prime ideals of ring R/B. Two cases arise:

**Case (1):** P = Q. In this case N(R/B) = P/B. It is clear now that if c+B is an element of c(N(R/B)), then c+B is a regular element of R/B. Hence by Goldie [2, Theorem 13], R/B has an artinian quotient ring.

**Case (2):** Q ≠ P. In this case observe that Q/B and P/B are distinct minimal prime ideals of the ring R/B. It is not difficult to see that in this case also that R/B has an artinian quotient ring.

**Theorem 9**. Let R be a Noetherian, link k-symmetric ring. Let R[X] be the polynomial ring over R. Let P be a prime ideal of R. Let P′ = P[X] be the extended prime ideal of P in the ring R[X]. Let Q′ be a prime ideal of R[X] such that Q′ → P′ is a link of the prime ideal Q′ to P′. Then Q′= Q[X] where Q = Q′ ∩ R. Moreover Q → P is a link of prime ideal Q of R to the prime ideal P of R.

**Proof.** Let Q′ ∩ P′/B′ be the linking bimodule for the link Q′ → P′, where B′ is an ideal of R[X]. Let B = B′ ∩ R. By hypothesis Q = Q′ ∩R. Clearly Q is a prime ideal of R whereas B is an ideal of R. Observe that by Lemma 1 any element d in c(Q) implies that d in c(Q′) also. Also it is obvious that any element c in c(P) implies that c in c(P′). We now have two cases.

**Case (1)** B < Q ∩ P. In this case by Lemma 2 Q ∩ P/B is a R/Q − R/P is a torsionfree bimodule. So Q → P is a right link between prime ideals Q and P of R. Again by Lemma 4 this induces a link Q[X] → P[X] of the prime ideals Q[X] and P [X] of R[X] via the ideal B[X]. Moreover R is right k-symmetric ring implies that $|R/P| = |R/Q|$. Hence by [3,Theorem13.17] $|R[X]/P[X]| = |R[X]/Q[X]|$. Now by theorem (8) we get that R[X]/B[X] has an Artinian quotient ring and since N (R[X]/B[X]) is a right k-homogeneous semiprime ideal of the ring R[X]/B[X], so by [1,Proposition 3.1], we get that N (R[X]/B[X]) is a right weakly ideal invariant ideal of R[X]/B[X]. Suppose $Q′ > Q[X]$. Then since N(R[X]/B[X]) is right weakly ideal invariant ideal of the ring R[X]/B[X], therefore $|Q′ ∩ P ′/Q′P ′| < |R[X]/P ′|$ because $|R[X]/Q′|<|R[X]/P [X]|$. This would immediately imply that $Q′ ∩ P ′/Q′P′$ is a right R[X]/P ′ torsion module. Hence in particular $Q′ ∩ P ′/B′$ is a right torsion R[X]/P ′ module. This is a contradiction to our hypothesis. Hence Q′ = Q[X] in this case.

**Case (2)** $B = Q ∩ P$. In this case $B[X] = Q[X] ∩ P [X]$ and hence B[X] is a semiprime ideal of the ring R[X]. Let S = set of regular elements of R[X]/B[X], then by Goldie [2, Theorem 13], we get that S is a left and a right ore-set of the ring R[X]/B[X]. Three subcases arise

**Subcase (a)** Suppose Q[X] and P [X] are distinct incomparable prime ideals. Then in this subcase it is clear that Q[X]/B[X] and P[X]/B[X] are the minimal prime ideals of the ring R[X]/B[X]. So S = c((P[X]) ∩ (Q[X])). Now the link $Q′ → P ′$ of the prime ideals Q′ and P ′ via the linking bimodule $Q′ ∩ P ′/B′$ can be thought of as a link in the ring R[X]/B[X] . But this is impossible by [3, Lemma 11.17].

**Subcase (b)** $Q[X]< P [X]$. In this subcase B[X] = Q[X]. Let $S_1$ = Set of regular elements of R/Q . Then it is not difficult to see that $S_1$ is a left and a right ore-set of R/B and hence also of the ring R[X]/B[X] (=R[X]/Q[X]). Now by Lemma 1, S1 is a subset of the set c((Q′)/B[X]). Now considering $Q′ → P ′$ as a link in the ring R[X]/B[X] we must have by [3,Lemma 12.17] that $S_1 = c((P ′)/B[X])$. But that is impossible because Q[X] < P [X]. Hence Q′ = Q[X].

The other subcase, that is when P [X] < Q[X], is on the same lines as subcase (b) above.
We must state that the proof is quite obvious in the case when Q[X]=P[X].

**Proposition 10**. Let R be a fully bounded Noetherian ring. Then R is a link k-symmetric ring.

**Proof**. See [3, Theorem 13.15].

**Theorem 11**. Let R be a fully bounded Noetherian ring. Let P be a fixed prime ideal of R. Then the following hold;
(a) If $Q′ → P [X]$ is a link of the prime ideal Q′ of R[X] to the extended prime ideal P [X] of R[X], then Q′ = Q[X], where $Q = Q′ ∩ R$.

(b) There are only finitely many prime ideals of R[X] right linked to P[X].
 Proof. (a) The proof of (a) follows by using Proposition 11 and Theorem 9.
(b) This follows immediately by [3, Theorem 13.22].

**Theorem 12**. Let R be a link k-symmetric Noetherian ring and let R[X] be the polynomial ring. Let $P_1$ and $Q_1$ be prime ideals of R[X] and
$Q_1 → P_1$ be a link of these prime ideals. Let $Q = Q_1 ∩ R$ and $P = P_1 ∩ R$. Then the following hold:

(a) If $P_1 = P [X]$, then $Q_1 = Q[X]$ and moreover $Q → P$ is a link of the
prime ideals Q and P of R.

(b) If $P_1 > P [X]$, then $Q_1 > Q[X]$.

Proof. (a) Proof of (a) follows directly from Theorems 7 and 9.
 (b) Proof of (b) is immediate from (a) above.

Remark:- We have come to know that our Theorem 9 is similar to Theorem 3.15 of the paper "Projective Prime ideals and localization in P.I rings" by A.W.chatters, C.R. Hajarnavis and R.M.Lissaman which appeared in J.London Math. Soc.(2)64 (2001) ( 1-12),L.M.S-2001 and in which the authors prove their result in the case when the base ring R is a P.I ring .In this respect we mention that our Theorem 9 is more general than the above mentioned Theorem 3.15. In fact the base ring R of our Theorem 9 which is a link krull symmetric noetherian ring generalizes noetherian P.I.rings and it is a significant open question whether or not a two sided noetherian ring is a link krull symmetric ring.